\documentclass[12pt]{article}
\input isolatin1.sty
 \usepackage{a4}
 \usepackage{amsmath}
 \usepackage{amssymb}
\usepackage{amsthm}
\usepackage{graphicx}
\newtheorem{thm}{Theorem}[section]
\newtheorem{cor}[thm]{Corollary}
\newtheorem{lem}[thm]{Lemma}
\newtheorem{prop}[thm]{Proposition}
\newtheorem{defn}[thm]{Definition}
\newtheorem{rem}[thm]{Remark}

\numberwithin{equation}{section} %%%%%%numera las ecuaciones (1.1), en
\newdimen\vgrowdim
\def\vgrow#1{\setbox0=\hbox{#1}\vgrowdim=\ht0 \advance\vgrowdim5pt
  \ht0=\vgrowdim\relax\vgrowdim=\dp0 \advance\vgrowdim5pt
  \dp0=\vgrowdim\box0}

\def\boxthem#1{\vbox{\def\\{\crcr\noalign{\hrule}}\offinterlineskip
  \halign{&\vrule##&\quad\hfil\vgrow{$##$}\quad\hfil\crcr#1\crcr}}}

\def\boxmatrix#1\endboxmatrix{\boxthem{#1}}

\title{The fixed point for a transformation of Hausdorff moment sequences
and iteration of a rational function
\footnote{This work was initiated during the visit of the second author to the University
of Copenhagen  partially supported by  D.G.E.S,
  ref. BFM2003-6335-C03-01, FQM-262 ({\em Junta de Andaluc\'{\i}a}), and by
  grant
  21-03-0324 from the Danish Natural Science Research Council.}}
\author{Christian Berg $^{\dagger}$ and Antonio J. Dur\'{a}n  $^{\ddagger}$\\
   \footnotesize $\dagger$ \ Institut for Matematiske Fag. K\o benhavns
   Universitet \\
    \footnotesize Universitetsparken 5; DK-2100 K\o benhavn \o,
   Denmark. berg@math.ku.dk \\  \footnotesize
   $\ddagger$ \  Departamento de An\'{a}lisis Matem\'{a}tico.
   Universidad de Sevilla \\
   \footnotesize Apdo (P. O. BOX) 1160. 41080 Sevilla. Spain. duran@us.es \\
\ \ }

 \date{\today}
\begin{document}
\maketitle

\begin{abstract}
We study the fixed point for a non-linear transformation in the set of
 Hausdorff moment sequences, defined by the formula:
$T((a_n))_n=1/(a_0+\cdots +a_n)$. We determine the corresponding measure
$\mu$, which has an increasing and convex density on $]0,1[$,
and we study some analytic functions related to it. The Mellin
transform  $F$ of $\mu$ extends to a meromorphic function
in the whole complex plane. It can be characterized in analogy 
with the Gamma function as the unique log-convex
function on $]-1,\infty[$ satisfying $F(0)=1$ and the functional
equation $1/F(s)=1/F(s+1)-F(s+1), s>-1$.
\end{abstract}

2000 {\it Mathematics Subject Classification}:
primary 44A60; secondary 30D05.

Keywords: Hausdorff moment sequence, iteration of rational functions

\section{Introduction and main results}
 Hausdorff moments sequences are sequences
of the form $\int _0^1t^nd \nu (t)$, $n\ge 0$, where $\nu $ is a
positive measure on $[0,1]$. 
Hausdorff moment sequences were characterized as completely monotonic sequences
in a fundamental paper by Hausdorff, see \cite{ha}.
For a recent study of Hausdorff moment sequences see \cite{df1},\cite{df2}.
Hausdorff moment sequences can also be characterized as bounded
Stieltjes moment sequences, where  Stieltjes moment sequences are 
of the form
$\int _0^{\infty} t^nd \nu (t)$, $n\ge 0$ for a
positive measure $\nu$ on $[0,\infty[$. 
 For a treatment of these concepts and the more general Hamburger moment
problem see the monograph by Akhiezer \cite{ak}.

In \cite{bd} the authors introduced  a non-linear multiplicative
 transformation from Hausdorff moment sequences to Stieltjes moment
 sequences.
  In \cite{bd2} we introduced a non-linear transformation
$T$ of the set of Hausdorff moment sequences into itself
by the formula:
\begin{equation}\label{eq:tms}
T((a_n))_n=1/(a_0+a_1+\cdots +a_n),\quad n\ge 0.
\end{equation}
The corresponding transformation of  positive measures on $[0,1]$ is 
denoted $\widehat{T}$. We recall from \cite{bd2} that if $\nu\neq 0$, then
$\widehat{T}(\nu)(\{0\})=0$ and
\begin{equation}\label{eq:bd2}
\int_0^1\frac{1-t^{z+1}}{1-t}\,d\nu(t)\int_0^1 t^z\,d\widehat{T}(\nu)(t)=1
\;\text{for}\; \Re z \geq 0.
\end{equation}
Assuming $\Re z>0$
we can consider $t^z=\exp(z\log t)$ as a continuous function on
$[0,1]$ with value 0 for $t=0$. Likewise
$(1-t^z)/(1-t)$ is
a continuous function for $t\in [0,1]$ with value $z$ for $t=1$.
If $\Re z=0,z\neq 0$ the function $t^z$ is only considered for $t>0$,
so it is important that $\widehat{T}(\nu)$ has no mass at zero. Finally
$t^0\equiv 1$.
It is clear that if $\nu $ is a probability measure, then
so is $\widehat{T}(\nu )$, and in this way we get a transformation of the
convex set of normalized Hausdorff moment sequences (i.e. $a_0=1$) as well as a
transformation of the set of probability measures on $[0,1]$. 
By Kakutani's theorem the transformation has a fixed point, and by
(\ref{eq:tms}) it is clear that a fixed point $(m_n)_n$ is uniquely
 determined by the
recursive equation
\begin{equation}\label{eq:fix}
m_0=1,\quad (1+m_1+\cdots +m_n)m_n=1,\quad n\ge 1.
\end{equation}                                    
Therefore
\begin{equation}\label{eq:fixrec}
m_{n+1}^2+\frac{m_{n+1}}{m_n}-1=0,
\end{equation}
giving
$$
m_1=\frac{-1+\sqrt 5}{2},\quad m_2=\frac{\sqrt{22+2\sqrt 5}-\sqrt
5-1}{4}, \cdots\,.
$$

The purpose of this paper is to study the Hausdorff moment sequence
$(m_n)_n$
and to determine its associated probability measure $\mu $, called the
{\it fixed point measure}.

We already know that $\mu(\{0\})=0$ because $\mu=\widehat{T}(\mu)$,
but it is also convenient to notice that $\mu(\{1\})=0$.
It is clear that $(m_n)_n$ decreases to $c=\mu(\{1\})\geq 0$, hence
$m_0+m_1+\ldots +m_n\geq (n+1)m_n$. By (\ref{eq:fix}) we get
$1\geq (n+1)m_n^2\geq (n+1)c^2$, showing that $c=0$.

In Section 4 we prove much more, namely
\begin{equation}\label{eq:sim}
m_n\sim 1/\sqrt{2n}\,\, \text{for}\,\, n\to\infty.
\end{equation}

 We will study $\mu$ by determining what we call the
 {\it Bernstein transform}
\begin{equation}\label{eq:Bern}
f(z)=\mathcal B(\mu)(z)=\int _0^1\frac{1-t^z}{1-t}d\mu (t),\quad  \Re z> 0
\end{equation}
as well as the {\it Mellin transform}
\begin{equation}\label{eq:Mellin}
F(z)=\mathcal M(\mu)(z)=\int _0^1t^zd\mu (t),\quad  \Re z> 0.
\end{equation}
These functions are clearly holomorphic in
the half-plane $\Re z>0$ and continuous in $\Re z\ge 0$, the latter because
$\mu(\{0\})=0$.

As a first result we prove:

\begin{thm}\label{thm:main-a} The functions $f,F$ can be extended to
 meromorphic functions in $\mathbb C$ and they satisfy
\begin{equation}\label{eq:fF=1}
f(z+1)F(z)=1,\quad z\in\mathbb C
\end{equation}
\begin{equation}\label{eq:z,z+1}
f(z)=f(z+1)-\frac{1}{f(z+1)},\quad z\in\mathbb C.
\end{equation}
They are holomorphic in $\Re z>-1$. Furthermore $z=-1$ is
a pole of $f$ and $F$.

The fixed point measure $\mu$ has the properties
\begin{equation}\label{eq:Landau}
\int_0^1 t^x\,d\mu(t)<\infty,\;\; x>-1;\quad 
\int_0^1 \frac{d\,\mu(t)}{t}=\infty.
\end{equation}
\end{thm}

\begin{proof} By (\ref{eq:bd2}) with $\nu$ replaced by the fixed point
measure $\mu$ we get $f(z+1)F(z)=1$ for $\Re z\ge 0$, showing in
particular that $f(z+1)$ and $F(z)$  are different from zero for $\Re z\ge 0$.
For $\Re z\ge 0$ we get by (\ref{eq:Bern}) 
$$
f(z+1)-f(z)=\int_0^1\frac{t^z-t^{z+1}}{1-t}\,d\mu(t)
=\int _0^1 t^{z}d\mu (t)=F(z)= \frac{1}{f(z+1)},
$$
which shows (\ref{eq:z,z+1}) for these values of $z$.

We remark that $\Re f(z)>0$ and in particular
$f(z)\neq 0$ for $\Re z>0$. This follows  by (\ref{eq:Bern}) because
$\Re(t^z)\le |t^z|<1$ for $0<t<1$ and $\Re z>0$. 

We next use  equation (\ref{eq:z,z+1}) to define $f(z)$ for $\Re z\ge
-1$, yielding a holomorphic continuation of $f$ to the open
half-plane $\Re z>-1$ because $f(z+1)\neq 0$.

Using equation (\ref{eq:z,z+1}) once more we obtain a meromorphic
extension of $f$ to the half-plane $\Re z>-2$. There will be poles at
points $z$ for which $f(z+1)=0$, in particular for $z=-1$ because
$f(0)=0$.

Repeated use of equation (\ref{eq:z,z+1}) makes it possible to obtain
a meromorphic extension to $\mathbb C$. At each step, $z$ will be a pole
if $z+1$ is a zero or a pole. 

At this stage we cannot give a complete picture of the poles of $f$,
but we return to that in Theorem \ref{thm:main-d}.

Having extended $f$ to a meromorphic function in $\mathbb C$ such that
(\ref{eq:z,z+1}) holds, we extend
$F$ to a meromorphic function in $\mathbb C$ such that equation
(\ref{eq:fF=1}) holds. 

Let us notice that also $F$ has no poles in $\Re z>-1$ 
because $f(z+1)\neq 0$. Moreover $z=-1$  is a pole of $F$ because $f(0)=0$.

By a classical result (going back to Landau for Dirichlet series), see
\cite[p. 58]{W}, we then get equation (\ref{eq:Landau}).
\end{proof}

The function $f$ can be characterized in analogy with the
Bohr-Mollerup theorem about the Gamma function, cf. \cite{ar}. More
precisely we have:

\begin{thm}\label{thm:main-b} The Bernstein transform (\ref{eq:Bern})
of the fixed point measure is a function  $f:\left]0,\infty\right[\to \left]0,\infty\right[$
 with the following properties
\begin{enumerate}
\item[(i)] $f(1)=1$,
\item[(ii)] $\log(1/f)$ is convex,
\item[(iii)] $f(s)=f(s+1)-1/f(s+1),\quad s>0$.
\end{enumerate}
Conversely, if $\tilde f:]0,\infty[\to ]0,\infty[$ satisfies
(i)-(iii), then it is equal to $f$ and
for $0<s\le 1$ we have
\begin{equation}\label{eq:iterate}
\tilde f(s)=\lim_{n\to\infty}\psi^{\circ n}\left(\frac{1}{m_{n-1}}\left(
\frac{m_{n-1}}{m_n}\right)^s\right),
\end{equation}
where $\psi$ is the rational function $\psi(z)=z-1/z$. In particular 
(\ref{eq:iterate}) holds for $f$.
\end{thm}

Here and elsewhere we use the notation for composition of mappings
$\psi^{\circ 1}(z)=\psi(z), \psi^{\circ
  n}(z)=\psi(\psi^{\circ(n-1)}(z)),\;n\ge 2$.
Theorem \ref{thm:main-b} will be proved in Section 3.
Using the relation $f(s+1)F(s)=1$ it is clear that Theorem
\ref{thm:main-b} can be reformulated to a characterization of $F$:

\begin{thm}\label{thm:main-c} There exists one and only one function 
$F:\left]-1,\infty\right[\to \left]0,\infty\right[$ with the following
properties
\begin{enumerate}
\item[(i)] $F(0)=1$,
\item[(ii)] $F$ is log-convex,
\item[(iii)] $1/F(s)=1/F(s+1)-F(s+1),\quad s>-1$,
\end{enumerate}
namely $F$ is the Mellin transform
$$
F(s)=\int_0^1 t^s\,d\mu(t),\quad s>-1
$$
of the fixed point measure.
\end{thm}

Let $\mathcal H$ denote the set of normalized Hausdorff moment
sequences considered as a subset of $[0,1]^{\mathbb N_0}$ with the
product topology, $\mathbb N_0=\{0,1,\ldots\}$.
In Section 2 we prove that the fixed point $\mathbf{m}=(m_n)_n$ is
attractive in the sense that for each $\mathbf{a}=(a_n)_n\in\mathcal H$
the sequence of iterates $T^{\circ n}(\mathbf{a})$ converges to
$\mathbf{m}$ in $\mathcal H$. Focusing on probability measures we see that  
every probability measure $\tau$ on
$[0,1]$ belongs to the domain of attraction of the fixed point measure
$\mu$ in the sense that 
$\lim_{n\to\infty} \widehat{T}^{\circ n}(\tau)=\mu$ weakly. 
For $q\in\mathbb R$ we denote by $\delta_q$ the probability measure 
with mass 1 concentrated at the point $q$. By specializing the
 iteration using $\tau=\delta_0$
we prove the following result:

\begin{thm}\label{thm:main-d} Let $f$ and $F$ be the meromorphic
  functions in $\mathbb C$ extending (\ref{eq:Bern}) and
  (\ref{eq:Mellin}) respectively.
The zeros and poles of $f$ are all simple and are contained in
$]-\infty,0]$. The zeros of $f$ are denoted $\xi_0=0$
 and $\xi_{p,k}, p=1,2,\ldots,
k=1,\ldots,2^{p-1}$ with $-p-1<\xi_{p,1}<\xi_{p,2}<\cdots
<\xi_{p,2^{p-1}}<-p$.

The poles of $f$ are $-l,\xi_{p,k}-l, l=1,2,\ldots$ with $p,k$ as
above.

Defining
\begin{equation}\label{eq:rho-f}
\rho _0=\frac{1}{f'(0)};\quad \rho _{p,k}=\frac{1}{f'(\xi _{p,k})},
\end{equation}
then $\rho_0,\rho_{p,k}>0$.

The following representations hold
\begin{equation}\label{eq:F}
F(z)=\frac{\rho_0}{z+1}+\sum_{p=1}^\infty\sum_{k=1}^{2^{p-1}}\frac{\rho_{p,k}}
{z+1-\xi_{p,k}},
\end{equation}
and
 \begin{equation}\label{eq:f}
f(z)=z\sum_{l=1}^\infty\left[\frac{\rho_0}{l(z+l)}
+\sum_{p=1}^\infty\sum_{k=1}^{2^{p-1}}\frac{\rho_{p,k}}
{(l-\xi_{p,k})(z+l-\xi_{p,k})}\right].
\end{equation}
The fixed point measure $\mu$ has an increasing and
convex density $\mathcal D$ with respect to Lebesgue measure on
$]0,1[$ and it is given by
\begin{equation}\label{eq:D}
\mathcal D(t)=\rho_0+ \sum_{p=1}^\infty\sum_{k=1}^{2^{p-1}}
\rho_{p,k} t^{-\xi_{p,k}}.
\end{equation}
\end{thm}

While clearly $\mathcal D(0)=\rho_0$, we prove in Theorem \ref{thm:one} that
$$
\mathcal D(t)\sim 1/\sqrt{2\pi(1-t)},\quad t\to 1.
$$
It is possible to obtain expressions for $\xi_{p,k}$ and $\rho_{p,k}$
in terms of the moments $(m_n)$. This is quite technical and is
given in Theorem \ref{thm:xi-rho}. 

We recall that a function $\varphi$ is called a {\it Stieltjes transform} if it can
be written in the form
\begin{equation}\label{eq:Sti}
\varphi(z)=a+\int_0^\infty\frac{d\sigma(x)}{x+z},\quad z\in\mathbb C\setminus
 ]-\infty,0],
 \end{equation}
where $a\ge 0$ and $\sigma$ is a positive measure on $[0,\infty[$ such that
(\ref{eq:Sti}) makes sense, i.e. $\int 1/(x+1)\,d\sigma(x)<\infty$.

It is clear that if $\sigma\neq 0$ then $\varphi$ is strictly decreasing
on $]0,\infty[$ with $a=\lim_{s\to\infty}\varphi(s)$. Furthermore, $\varphi$
is holomorphic in $\mathbb C\setminus ]-\infty,0]$ with
$$
\frac{\Im \varphi(z)}{\Im z} < 0\;\;\text{for}\;\; z \in \mathbb
C\setminus\mathbb R,
$$
so in particular $\varphi$ is never zero in $\mathbb C\setminus ]-\infty,0]$.
The Stieltjes transforms we are going to consider will be meromorphic in
$\mathbb C$. The function (\ref{eq:Sti}) is meromorphic
precisely when the measure $\sigma$ 
  is discrete and the set of mass-points have no finite 
 accumulation points, i.e. if and only if
 $$
 \varphi(z)=a + \sum_{p=0}^\infty\frac{\sigma_p}{z+\eta_p}
 $$
 with $\sigma_p>0,\;0\le \eta_0<\eta_1<\eta_2<\ldots\to\infty$.

For results about Stieltjes transforms see \cite{bf}. Stieltjes transforms
are closely related to Pick functions, cf. \cite{ak},\cite{d}.
We recall that a Pick function is a holomorphic function $\varphi:\mathbb
C\setminus\mathbb R\to\mathbb C$ satisfying
$$
\frac{\Im \varphi(z)}{\Im z} \ge 0\;\;\text{for}\;\; 
z \in \mathbb C\setminus\mathbb R,
$$
so if $\varphi\ne 0$ is a Stieltjes transform, then $1/\varphi$ is a
Pick function. Notice that $z/(z+a)$ is a Pick function for any $a>0$.

\begin{cor}\label{thm:BP} In the notation of Theorem \ref{thm:main-d}
$f(z)/z$ and $F(z)$  are Stieltjes transforms and $f$ is a Pick function.
\end{cor}

We have used the name Bernstein transform for (\ref{eq:Bern}). In
general, if $\nu$ is a positive finite measure on
$]0,1]$, we call
\begin{equation}\label{eq:Berngen}
\mathcal B(\nu)(z)=\int _0^1\frac{1-t^z}{1-t}\,d\nu(t)
\end{equation}
the Bernstein transform of $\nu$, because it
is a Bernstein function in the terminology of \cite{bf}. In fact we
can write
$$
\mathcal B(\nu)(z)=\nu(\{1\})z+
\int _0^\infty \left(1-e^{-xz}\right) d\lambda (x), \quad \Re z\ge 0,
$$
where $\lambda$ is  defined as the image measure
of $(1-t)^{-1}(\nu|]0,1[)$ under $\log(1/x)$ mapping $]0,1[$ onto
$]0,\infty[$. We recall that $\lambda$ is called the L\'evy measure
of the Bernstein function. It follows that $\mathcal B(\nu)'$ is a
completely monotonic function.
Bernstein functions are very important in the theory of L\'evy
processes, see \cite{bt}.

In section 4 we prove that $(m_n)_n$ is infinitely divisible in the
sense that $(m_n^\alpha)_n$ is a Hausdorff moment sequence for all
 $\alpha>0$. 

\section{An iteration leading to the fixed point measure}

For $n=0,1,\ldots$ we denote the moments of
 $\mu_n=\widehat{T}^{\circ n}(\delta_0)$ by
$(m_{n,k})_k$, i.e.
$$
\int_0^1 t^k\,d\widehat{T}^{\circ n}(\delta_0)(t)=m_{n,k},
$$
hence for $n\ge 1$
\begin{equation}\label{eq:m}
m_{n,k}=\left(m_{n-1,0}+m_{n-1,1}+\cdots +m_{n-1,k}\right)^{-1}.
\end{equation}
 Notice that $m_{n,0}=1$ for all $n$ and $m_{0,k}=\delta_{0k},
 m_{1,k}=1, m_{2,k}=1/(k+1)$ for all $k$.

\begin{lem}\label{thm:limit} For fixed $k=0,1,\ldots$ we have
$$
m_{0,k}\leq m_{2,k}\leq m_{4,k}\leq\ldots
$$
$$
m_{1,k}\geq m_{3,k}\geq m_{5,k}\geq\ldots
$$
and these sequences have the same limit
$$
\lim_{n\to\infty} m_{2n,k}=\lim_{n\to\infty} m_{2n+1,k}=m_k,
$$
where $(m_k)_k$ is the fixed point given by (\ref{eq:fix}). 

 Furthermore, $\lim_{k\to\infty} m_{n,k}=0$ for $n\ge 2$, implying that       
 $\mu_n=\widehat{T}^{\circ n}(\delta_0)$ has no mass at $t=1$ for $n\ge 2$.
\end{lem}

\begin{proof} Since the result is trivial for $k=0$, we assume that
  $k\geq 1$ and have
$$
0=m_{0,k}< m_{2,k}=\frac{1}{k+1};\quad 1=m_{1,k}>m_{3,k}
=\frac{1}{\mathcal H_{k+1}},
$$
where $\mathcal H_p=1+\tfrac12+\cdots+\tfrac1p$ is the $p$'th harmonic
number. We now get
$$
\frac{1}{m_{4,k}}=\sum_{j=0}^k m_{3,j}< k+1
$$
hence $m_{4,k}>m_{2,k}$. We next use this to conclude
$$
\frac{1}{m_{5,k}}=\sum_{j=0}^k m_{4,j} >\sum_{j=0}^k
m_{2,j}=\frac{1}{m_{3,k}},
$$
hence $m_{5,k}<m_{3,k}$. It is clear that this procedure can be
continued and reformulated to a proof by induction.

Defining
$$
m_k'=\lim_{n\to\infty} m_{2n,k},\quad m_k''=\lim_{n\to\infty} m_{2n+1,k},
$$
we get the following relations from (\ref{eq:m})
\begin{equation}\label{eq:lim}
m_k'=\left(1+m_1''+\cdots +m_k''\right)^{-1},\quad
m_k''=\left(1+m_1'+\cdots +m_k'\right)^{-1},\quad k\geq 1,
\end{equation}
because clearly $m_0'=m_0''=m_0=1$. It follows easily by induction
using (\ref{eq:lim}) that
$m_k'=m_k''=m_k$ for all $k$.

 Since $m_{2n,k}\le m_k$
we get $\lim_{k\to\infty}m_{2n,k}=0$. Furthermore, for $n\ge 1$
$$
\frac{1}{m_{2n+1,k}}=\sum_{j=0}^k m_{2n,j}\ge \sum_{j=0}^k m_{2,j}=
\mathcal H_{k+1}
$$
and hence $\lim_{k\to\infty} m_{2n+1,k}=0.$
\end{proof}

We recall that $\mathcal H$ denotes the set of normalized Hausdorff moment
sequences $\mathbf{a}=(a_n)_n$. The mapping $\nu\to (\int
x^n\,d\nu(x))_n$ from the set of probability measures $\nu$ on $[0,1]$
to $\mathcal H$ is a homeomorphism between compact sets, when the set
 of probability
measures carries the weak topology and $\mathcal H$ carries the
topology inherited from $[0,1]^{\mathbb N_0}$ equipped with the
product topology.

Defining an order relation $\le$ on $\mathcal H$ by
writing $\mathbf{a}\le\mathbf{b}$ if $a_k\le b_k$ for $k=0,1,\ldots$,
we easily get the following Lemma:

\begin{lem}\label{thm:Tab}
The transformation $T:\mathcal H\to\mathcal H$ is decreasing,
i.e. 
$$
\mathbf{a}\le \mathbf{b}\Rightarrow T(\mathbf{a})\ge T(\mathbf{b}).
$$
\end{lem}

\begin{thm}\label{thm:basin} For every $\mathbf{a}\in\mathcal H$ we have
$$
\lim_{n\to\infty}T^{\circ n}(\mathbf{a})=\mathbf{m},
$$
where $\mathbf{m}=(m_n)_n$ is the fixed point.
\end{thm}

\begin{proof} For $0\le q\le 1$ we write
  $\pmb{\underline{q}}=(q^n)_n$, hence $\pmb{\underline{0}}\le\mathbf{a}
\le \pmb{\underline{1}}$ for every $\mathbf{a}\in\mathcal H$. By
 Lemma \ref{thm:Tab} we get
$$
T^{\circ (2n)}(\pmb{\underline{0}})\le T^{\circ (2n)}(\mathbf{a})
\le T^{\circ (2n)}(\pmb{\underline{1}})=T^{\circ
  (2n+1)}(\pmb{\underline{0}})
$$
$$
T^{\circ (2n+1)}(\pmb{\underline{0}})\ge T^{\circ (2n+1)}(\mathbf{a})
\ge T^{\circ (2n+1)}(\pmb{\underline{1}})=T^{\circ (2n+2)}(\pmb{\underline{0}}),
$$
and since $\lim_{n\to\infty} T^{\circ
  n}(\pmb{\underline{0}})=\mathbf{m}$ by Lemma \ref{thm:limit}, we get
$$
 \lim_{n\to\infty}T^{\circ (2n)}(\mathbf{a})= \lim_{n\to\infty}T^{\circ
   (2n+1)}(\mathbf{a})=\mathbf{m}.
$$
\end{proof} 
 
Theorem \ref{thm:basin} can also be expressed that $\widehat{T}^{\circ
  n}(\tau)\to\mu$ weakly for any probability measure $\tau$ on $[0,1]$.
Specializing this to $\tau=\delta_0$ and using  formula (\ref{eq:bd2}),
we obtain:

\begin{cor}\label{thm:cor1} The iterated sequence $\mu_n=\widehat{T}^{\circ
    n}(\delta_0)$ of measures converges weakly to the fixed point
  measure $\mu$ and
\begin{equation}\label{eq:mun}
\int_0^1\frac{1-t^{z+1}}{1-t}\,d\mu_n(t)\int_0^1
t^z\,d\mu_{n+1}(t)=1,\quad \Re z\geq 0,n=0,1,\ldots
\end{equation}
\end{cor}

We have $\mu_0=\delta_0, \mu_1=\delta_1, \mu_2=\chi_{]0,1[}(t)dt$,
where $\chi_{]0,1[}(t)$ denotes the indicator function for the interval
$]0,1[$. The Bernstein transform of the measure $\mu_2$ is
\begin{equation}\label{eq:B2}
\mathcal B(\mu_2)(z)=\int_0^1\frac{1-t^z}{1-t}\,dt=
\sum_{l=1}^\infty \left(\frac{1}{l}-\frac{1}{z+l}\right)=
\Psi(z+1)+\gamma,
\end{equation}
where $\gamma$ is Euler's constant and
$\Psi(x)=\Gamma'(x)/\Gamma(x)$ is the Digamma function.

 The measure $\mu_3$ has been calculated in \cite{bd2} and
 the result is
$$
\mu_3=\left(\sum_{p=0}^\infty \alpha_p
  t^{-\xi_p}\right)\chi_{]0,1[}(t)dt,
$$
where $0=\xi_0>\xi_1>\xi_2>\ldots$ satisfy $-p-1<\xi_p<-p$ for
$p=1,2,\ldots$ and $\alpha_p>0, p=0,1,\ldots.$ More precisely, it was
proved that $\xi_p$ is the unique solution $x\in \left]-p-1,-p\right[$ of the
equation $\Psi(1+x)=-\gamma$. Writing
$\xi_p=-p-1+\delta_p$, we have $0<\delta_{p+1}<\delta_p<\tfrac12,\;\delta_p
\sim 1/\log p,\; p\to\infty$. Furthermore,
$\alpha_p=1/\Psi'(1+\xi_p)\sim 1/\log^2 p$. Since $\sum
\alpha_p/(1-\xi_p)=1$, we have the crude estimate $\alpha_p<p+2$.

We shall now prove that all the measures $\mu_n,n\geq 4$ have a
form similar to that of $\mu_3$.

\begin{lem}\label{thm:munmun}
For $n\geq 3$ the measure $\mu_n$ has the form
\begin{equation}\label{eq:mun1}
\mu_n=\left(\rho_0^{(n)}+\sum_{p=1}^\infty\sum_{k=1}^{N(n,p)}
\rho_{p,k}^{(n)} t^{-\xi_{p,k}^{(n)}}\right)\chi_{]0,1[}(t)dt,
\end{equation}
where for each $p\geq 1$ 
\begin{enumerate}
\item[(i)] $1\leq N(n,p)\leq 2^{p-1}$,
\item[(ii) ]$-p-1<\xi_{p,1}^{(n)}<\xi_{p,2}^{(n)}< \ldots
<\xi_{p,N(n,p)}^{(n)}<-p$,
\item[(iii)] $0<\rho_0^{(n)}<1$, $0<\rho_{p,k}^{(n)}<p+2,
\quad k=1,\ldots,N(n,p)$.
\end{enumerate}
\end{lem}

\begin{proof}
The result for $n=3$ follows from the description above from
\cite{bd2} with $\rho_0^{(3)}=\alpha_0, N(3,p)=1,
\rho_{p,1}^{(3)}=\alpha_p, \xi_{p,1}^{(3)}=\xi_p$.

Assume now that the result holds for $\mu_n$ and let us prove it
 for $\mu_{n+1}$. For $\Re z>0$ we then have

$$
f_n(z):=\mathcal B(\mu_n)(z)=\int _0^1\frac{1-t^z}{1-t}\,d\mu_n(t)=
\sum_{l=0}^\infty\int_0^1\left(t^l-t^{z+l}\right)\,d\mu_n(t)
$$
$$
=\sum_{l=0}^\infty\left[\rho_0^{(n)}\int_0^1\left(t^l-t^{z+l}\right)\,dt+
\sum_{p=1}^\infty\sum_{k=1}^{N(n,p)}\rho_{p,k}^{(n)}
\int_0^1\left(t^{l-\xi_{p,k}^{(n)}}-t^{z+l-\xi_{p,k}^{(n)}}\right)\,dt\right]
$$
$$
=z\sum_{l=1}^\infty\left[\frac{\rho_0^{(n)}}{l(z+l)}
+\sum_{p=1}^\infty\sum_{k=1}^{N(n,p)}\frac{\rho_{p,k}^{(n)}}
{(l-\xi_{p,k}^{(n)})(z+l-\xi_{p,k}^{(n)})}\right].
$$

This shows that $f_n(z)/z$ is a Stieltjes transform and a meromorphic
function in $\mathbb C$ with poles at the points
$$
-l, \xi_{p,k}^{(n)}-l,\quad l=1,2,\ldots, p=1,2,\ldots,
k=1,\ldots,N(n,p),
$$
so in the interval $]-p-1,-p]$ we have the poles
\begin{equation}\label{eq:N(n,p)}
-p,\;\xi_{p-l,k}^{(n)}-l, k=1,\ldots,N(n,p-l),l=1,\ldots,p-1.
\end{equation}
 Since $f_n(x)/x$ is strictly decreasing
between the poles, we conclude that there is precisely one simple zero
between two consecutive poles.
Let $N(n+1,p)$ denote the number of zeros of $f_n$ in $]-p-1,-p[$ and let
$\xi_{p,k}^{(n+1)}$ denote the zeros numbered such that
$$
-p-1<\xi_{p,1}^{(n+1)}<\xi_{p,2}^{(n+1)}<\ldots<\xi_{p,N(n+1,p)}^{(n+1)}<-p.
$$
In addition also $z=0$ is a zero of $f_n$. There are no  zeros
or poles in $\mathbb C\setminus ]-\infty,0]$ because $f_n(z)/z$ is a
 Stieltjes transform.

We are now ready to prove equation (\ref{eq:mun1}) and (i)--(iii)
 with $n$ replaced by $n+1$.

(i). By (\ref{eq:N(n,p)}) we get
$$
N(n+1,p)\le  1+\sum_{l=1}^{p-1} N(n,p-l)
\leq 1+\sum_{l=1}^{p-1} 2^{p-l-1}=2^{p-1}.
$$

(ii) is clear by definition, when we have proved that the measure
$\mu_{n+1}$ has the form (\ref{eq:mun1}) using the numbers 
$\xi^{(n+1)}_{p,k}$.

(iii).  By a
 classical result, see \cite{r},\cite{i},\cite{b},  $1/f_n(z)$ is a
Stieltjes transform because $f_n(z)/z$ is so, i.e.
$$
\frac{1}{f_n(z)}=\frac{\rho_0^{(n+1)}}{z}+
\sum_{p=1}^\infty\sum_{k=1}^{N(n+1,p)}\frac{\rho_{p,k}^{(n+1)}}
{z-\xi_{p,k}^{(n+1)}},
$$
with $\rho_0^{(n+1)},\rho_{p,k}^{(n+1)}>0$.
There is no constant term in the Stieltjes representation because
 $f_n(x)\to\infty$ for $x\to\infty$. In fact, by Lemma \ref{thm:limit} we get
 $$
 \lim_{x\to\infty} f_n(x)=\int_0^1\frac{d\mu_n(t)}{1-t}=\sum_{k=0}^\infty
 m_{n,k}=\lim_{k\to\infty}\frac{1}{m_{n+1,k}}=\infty.
 $$
Note that
\begin{equation}\label{eq:rho-n}
\rho_0^{(n+1)}=\frac{1}{f_n'(0)},\quad \rho_{p,k}^{(n+1)}=
\frac{1}{f_n'(\xi_{p,k}^{(n+1)})}.
\end{equation}

By (\ref{eq:mun}) we get 
$$
\int_0^1 t^z\,d\mu_{n+1}(t)=\frac{1}{f_n(z+1)}=
\frac{\rho_0^{(n+1)}}{z+1}+
\sum_{p=1}^\infty\sum_{k=1}^{N(n+1,p)}\frac{\rho_{p,k}^{(n+1)}}
{z+1-\xi_{p,k}^{(n+1)}},
$$
which shows that
$$
\mu_{n+1}=\left(\rho_0^{(n+1)}+\sum_{p=1}^\infty\sum_{k=1}^{N(n+1,p)}
\rho_{p,k}^{(n+1)} t^{-\xi_{p,k}^{(n+1)}}\right)\chi_{]0,1[}(t)dt,
$$
which is (\ref{eq:mun1}) with $n$ replaced by $n+1$.

Since $\mu_{n+1}$ is a probability measure we get 
$$
\rho_0^{(n+1)}<1,\quad \rho_{p,k}^{(n+1)}\int_0^1
t^{-\xi_{p,k}^{(n+1)}}\,dt<1,
$$
hence
$$
\rho_{p,k}^{(n+1)} < 1-\xi_{p,k}^{(n+1)}<p+2.
$$ 
\end{proof}

\begin{cor}\label{thm:cor2} For $n\ge 0$ let
 $\mu_n=\widehat{T}^{\circ n}(\delta_0)$.
The functions $f_n=\mathcal B(\mu_n)$ are
meromorphic Pick functions and the functions $F_n=\mathcal M(\mu_n)$
are meromorphic Stieltjes transforms satisfying
\begin{equation}\label{eq:n-n+1}
f_n(z+1)F_{n+1}(z)=1,\quad z\in\mathbb C.
\end{equation}
All zeros and poles of $f_n$ are contained in $]-\infty,0]$.
\end{cor}

\begin{proof} We have $f_0(z)=1, f_1(z)=z, F_0(z)=0, F_1(z)=1,
F_2(z)=1/(z+1)$
and for $n\ge 2$ the result follows from  
Lemma \ref{thm:munmun} and its proof.
\end{proof}

In order to obtain a limit result for $n\to\infty$ in 
Corollary \ref{thm:cor2} we need the following:

\begin{lem}\label{thm:Stieltjes} Let $(\varphi_n)_n$ be a sequence of Stieltjes
transforms of the form
$$
\varphi_n(z)=\int_0^\infty \frac{d\sigma_n(x)}{x+z},\quad n=1,2,\ldots
$$
and assume that $\varphi_n(z)\to\varphi(z)$ uniformly on compact subsets of
$\Re z> 0$ for some holomorphic function $\varphi$ on the right half-plane.

Then $\varphi$ is a Stieltjes transform
$$
\varphi(z)=a+\int_0^\infty \frac{d\sigma(x)}{x+z}
$$
and $\lim_{n\to\infty}\sigma_n=\sigma$ vaguely. Furthermore,
$\varphi_n(z)\to\varphi(z)$ uniformly on compact subsets of
$\mathbb C\setminus ]-\infty,0]$.
\end{lem}

\begin{proof} Since
$$
\int_0^\infty \frac{d\sigma_n(x)}{x+1}=\varphi_n(1)\to\varphi(1),
$$
there exists a constant  $K>0$ such that $\int 1/(x+1)\,d\sigma_n(x)\le K$
for all $n$. Let $\sigma$ be a vague accumulation point for $(\sigma_n)_n$.
Replacing $(\sigma_n)_n$ by a subsequence we can assume without loss of generality
that $\sigma_n\to\sigma$ vaguely. By standard results in measure theory,
cf. \cite[Prop. 4.4]{bcr}, we  have
$$
\int_0^\infty \frac{d\sigma(x)}{x+1}\le K,\quad\lim_{n\to\infty}
\int f\,d\sigma_n=\int f\,d\sigma
$$
for any continuous function $f:[0,\infty[\to\mathbb C$ which is $o(1/(x+1))$
for $x\to\infty$. In particular
$$
\varphi_n'(z)=-\int_0^\infty \frac{d\sigma_n(x)}{(x+z)^2}\to
-\int_0^\infty \frac{d\sigma(x)}{(x+z)^2},\quad z\in
 \mathbb C\setminus ]-\infty,0],
$$
showing that
$$
\varphi'(z)=-\int_0^\infty \frac{d\sigma(x)}{(x+z)^2},\quad \Re z>0,
$$
hence
$$
   \varphi(z)=a+\int_0^\infty \frac{d\sigma(x)}{x+z},\quad \Re z>0
$$                                                                
for some constant $a$.
Using $\varphi(x)=\lim_{n\to\infty}\varphi_n(x)\ge 0$ for $x>0$, we get
$a\ge 0$, showing that $\varphi$ is a Stieltjes transform. By unicity
of $a$ and $\sigma$
in the representation of $\varphi$ as a Stieltjes transform, we conclude
that the accumulation point $\sigma$ is unique, hence
 $\lim_{n\to\infty} \sigma_n=\sigma$ vaguely.

 It is now easy to see that $(\varphi_n(z))_n$ is uniformly bounded on compact
 subsets of  $\mathbb C\setminus ]-\infty,0]$, and the last assertion
 of Lemma \ref{thm:Stieltjes} is a
  consequence of the Stieltjes-Vitali theorem.
\end{proof}

\medskip
{\it Proof of Theorem \ref{thm:main-d}.}

\medskip   
 From Lemma \ref{thm:munmun} follows
  that the Mellin transform $\mathcal M(\mu_n)(z)$ coincides on $\Re z\ge 0$
  with the meromorphic function
$$
\frac{\rho_0^{(n)}}{z+1}+\sum_{p=1}^\infty\sum_{k=1}^{N(n,p)}
\frac{\rho_{p,k}^{(n)}}
{z+1-\xi_{p,k}^{(n)}}=\int_0^\infty \frac{d\sigma_n(x)}{x+z},
$$
where $\sigma_n$ is the discrete measure
$$
  \sigma_n=\rho_0^{(n)}\delta_1+\sum_{p=1}^\infty\sum_{k=1}^{N(n,p)}
\rho_{p,k}^{(n)} \delta_{1-\xi_{p,k}^{(n)}}.
$$

Since $\mathcal M(\mu_n)(z)\to\mathcal M(\mu)(z)$ uniformly on compact
 subsets of $\Re z>0$ by Corollary \ref{thm:cor1}, it follows by
 Lemma \ref{thm:Stieltjes} that $\mathcal M(\mu)$ is a Stieltjes transform
$$
\mathcal M(\mu)(z)=a+\int_0^\infty\frac{d\sigma(x)}{x+z},
$$
and $\sigma_n\to \sigma$ vaguely.
Since $\mathcal M(\mu)(k)=m_k\to 0$ as $k\to\infty$, 
 we get $a=0$.
Using that $\sigma_n$ has at most $2^{p-1}$
 mass points
in $[p+1,p+2],\,p=1,2,\ldots$ and that $\rho_{p,k}^{(n)}< p+2$ by
Lemma \ref{thm:munmun}, we can write
$$
  \sigma=\rho_0\delta_1+\sum_{p=1}^\infty\sum_{k=1}^{N_p}
\rho_{p,k} \delta_{1-\xi_{p,k}},
$$
with $\rho_0\ge 0,0<\rho_{p,k}\le p+2$ and
$-p-1\le\xi_{p,1}<\xi_{p,2}<\cdots<\xi_{p,N_p}<-p$, where $N_p\le 2^{p-1}$.
At this stage we cannot confirm that $\rho_0>0, -p-1<\xi_{p,1}$,
 $N_p=2^{p-1}$ and that $\xi_{p,k}$ are the zeros of $f$.
The function
\begin{equation}\label{eq:Fprel}
\frac{\rho_0}{z+1}+\sum_{p=1}^\infty\sum_{k=1}^{N_p}\frac{\rho_{p,k}}
{z+1-\xi_{p,k}}
\end{equation}
is a meromorphic extension of $\mathcal M(\mu)$  and therefore equal
to the meromorphic function $F$ of Theorem \ref{thm:main-a}. This shows that
$\mu$ has the density
\begin{equation}\label{eq:denshelp}
\mathcal D(t)=\rho_0+ \sum_{p=1}^\infty\sum_{k=1}^{N_p}
\rho_{p,k} t^{-\xi_{p,k}},
\end{equation}
which is clearly increasing and convex since $-\xi_{p,k}\ge -1$.
Finally, by (\ref{eq:denshelp}) the Bernstein transform $\mathcal
B(\mu)$  has the meromorphic extension
\begin {equation}\label{eq:fprel}
z\sum_{l=1}^\infty\left[\frac{\rho_0}{l(z+l)}
+\sum_{p=1}^\infty\sum_{k=1}^{N_p}\frac{\rho_{p,k}}
{(l-\xi_{p,k})(z+l-\xi_{p,k})}\right],
\end{equation}
which is a Pick function. The function given by (\ref{eq:fprel}) equals
the meromorphic function $f$ of Theorem \ref{thm:main-a}.
 By Lemma \ref{thm:Stieltjes} applied to the Stieltjes
 transforms $f_n(z)/z$, we conclude that $f_n(z)\to f(z)$ uniformly on
 compact subsets of $\mathbb C\setminus ]-\infty,0]$.

We already know from  Theorem \ref{thm:main-a} that $F$ has a pole at
 $z=-1$ and hence
$\rho_0>0$. The remaining poles of $F$ are $\xi_{p,k}-1$, so by
formula (\ref{eq:fF=1}) the zeros of $f$ are $z=0$ and $z=\xi_{p,k}$.
By the expression (\ref{eq:fprel}) for $f$
the  poles of $f$ are $-l,\xi_{p,k}-l$ and therefore
 $-p-1<\xi_{p,1},\,p=1,2,\ldots$. 

We have now proved that the zeros and poles of $f$ are all simple and
are contained in $]-\infty,0]$. Since $f(z+1)F(z)=1$ we get by
 (\ref{eq:Fprel}) that
$$
\frac{1}{f(z)} = \frac{\rho_0}{z}+\sum_{p=1}^\infty\sum_{k=1}^{N_p}
\frac{\rho_{p,k}}{z-\xi_{p,k}},
$$
which shows equation (\ref{eq:rho-f}).

To finish the proof we shall establish that $N_p=2^{p-1}$.

From the functional equation (\ref{eq:z,z+1}) and the fact that $f$ is strictly
increasing between the poles, we see the following about the
generation of zeros and poles of $f$:
\begin{enumerate}
\item If $z+1$ is regular point, then $f(z+1)=\pm 1$ if  and only if $f(z)=0$.
\item If  $z+1$ is regular point, then 
$f(z+1)=0$ if and only if $z$ is a pole. In the affirmative case
$\text{Res}(f,z)=-1/f'(z+1)$.
\item If $z+1$ is a pole then $z$ is a pole with the same residue as in $z+1$.
\item For a pole $\beta$ let $\alpha_\beta$ be the smallest zero in
  $]\beta,\infty[$. Then $f(]\beta,\alpha_\beta[)=]-\infty,0[$ and
there exists a unique point $x_*$ in $]\beta,\alpha_\beta[$ such that
$f(x_*)=-1$.
\item For a pole $\beta$ let $\gamma_\beta$ be the biggest zero in
 $]-\infty,\beta[$. Then $f(]\gamma_\beta,\beta[)=]0,\infty[$ and
 there exists a unique point $x^*$ in $]\gamma_\beta,\beta[$ such that 
$f(x^*)=1$. 
\end{enumerate}

From 1.-5. we  deduce that $f$ has
the following properties. Since $f(0)=0$ we see that $f$ has poles at
 $z=-1,-2,\ldots$ in accordance with (\ref{eq:fprel}). There are no poles
 in $]-2,-1[$ since $f$ is
 regular in $]-1,0[$ and non-zero. Notice that $f$ is strictly
 increasing on $\left]-1,\infty\right[$ mapping this interval onto the
 whole real line by (\ref{eq:fprel}).
There is a unique point $x_*\in\left]-1,0\right[$ such that
 $f(x_*)=-1$, hence $x_*-1$ is a zero and $x_*-2,x_*-3,\ldots$ are
 poles. In $]-3,-2]$ there are two poles namely $x_*-2$ and $-2$ and
 since  $f$ is strictly increasing between consecutive poles we have
 two zeros in $]-3,-2[$. By
 induction it is easy to see that there are exactly $2^{p-1}$
poles in each interval $]-p-1,-p]$ and $2^{p-1}$  zeros in the open interval
 $\left]-p-1,-p\right[$, $p\ge 1$. This shows that $N_p=2^{p-1}$.
Note that $\xi_{1,1}=x_*-1$.
\hfill$\square$

\begin{figure}[htb]
\begin{center}
\includegraphics[scale=0.5]{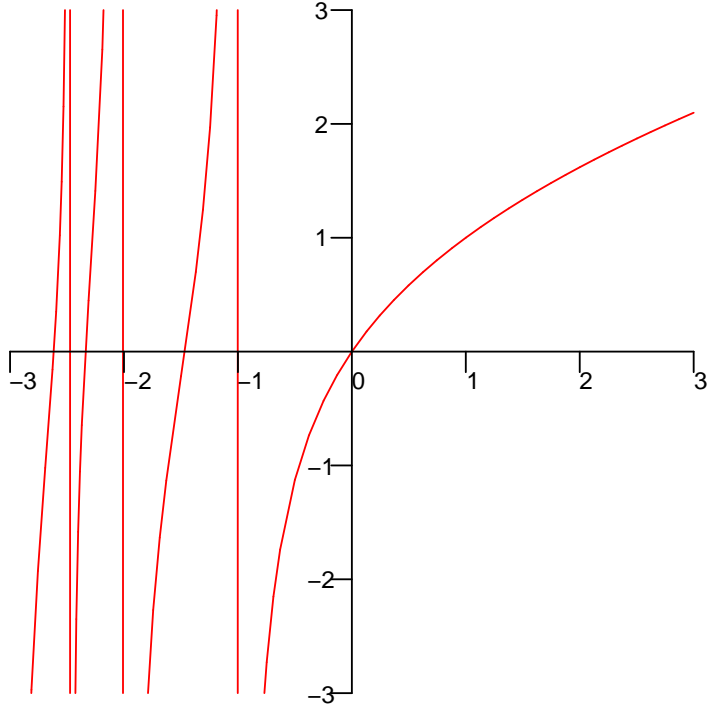}
\end{center}
\caption{The graph of $f$ with vertical lines at the poles}
\end{figure}
     
\medskip
We give some further information about the poles of $f$.

We call the negative integers poles of the {\it first generation}  of
$f$ and say
 that a pole of $f$ is of the
 {\it $l$-th generation}, $l\ge 2$, if it is generated by a zero 
$\xi _{l-1,k}$,
i.e. the pole is of the form $\xi _{l-1,k}-m$, for some integer
$m\ge 1$.
Then it can easily be proved by induction on $p$ that:

 \begin{enumerate}
\item In $]-p-1,-p]$ there is one pole of the first generation
(namely, $-p$), one pole of the
second generation
(namely $\xi _{1,1}-p+1$), and for $l=3,\ldots , p$, $2^{l-2}$ poles of
 the $l$-th generation (so that
the total number of poles is $1+\sum _{l=2}^{p}2^{l-2}=2^{p-1}$).
\item For each interval $[-p-1,-p]$,
the poles of one generation separate the set of poles of lower generations,
and the zeros $\xi _{p,k}$, $k=1,\ldots , 2^{p-1}$, separate the set of all poles.
That means that the set of poles of generation less than or equal to $l$
separate the zeros $\xi _{p,k}$, $k=1,\ldots , 2^{p-1}$, in groups of
$2^{p-l}$ consecutive elements.
\item For $l\geq 2$ the poles in $\left]-p-1,-p\right[$ of the $l$-th
 generation are zeros of $f(z+p-l+1)$ but they are still poles of $f(z+j)$ if
$0\le j\le p-l$.
\end{enumerate}

\section{Iteration of the rational
 function $\psi$}

In this section we will prove Theorem \ref{thm:main-b} and discuss the
relationship with the classical study of iteration of rational
functions of degree $\ge 2$, cf. e.g. \cite{bear}.

We have already introduced the rational function $\psi$ by
\begin{equation}\label{eq:psi}
\psi(z)=z-\frac{1}{z}.
\end{equation}
It is a mapping of $\mathbb C\setminus\{0\}$ onto $\mathbb
C$ with a simple pole at $z=0$. Moreover,
$\psi(0)=\psi(\infty)=\infty$. It is two-to-one with the exception
that $\psi(z)=\pm 2i$ has only one solution $z=\pm i$.
 It is strictly increasing on the
half-lines $\left]-\infty,0\right[$ and $\left]0,\infty\right[$,
mapping each of them onto $\mathbb R$.
The functional equation (\ref{eq:z,z+1}) can be written
\begin{equation}\label{eq:it}
f(z)=\psi(f(z+1)).
\end{equation}

We notice that $\psi$ and hence all iterates $\psi^{\circ n}$
are Pick functions. It is convenient to define $\psi^{\circ 0}(z)=z$.
We claim that the Julia set is
$J(\psi)=\mathbb R^*$,
and the Fatou set is $F(\psi)=\mathbb C\setminus\mathbb R$. This is
because $\psi$ is conjugate to the rational function 
$$
R(z)=\frac{3z^2+1}{z^2+3}
$$
i.e. $g\circ R=\psi\circ g$, where $g$ is the M\"obius transformation 
$g(z)=i(1+z)/(1-z)$. Note that $g$ is the Cayley transformation
mapping the unit circle $\mathbb T$ onto
$\mathbb R^*$. In \cite[p.200]{bear} the Julia set of $R$ is
determined as $J(R)=\mathbb T$, and the assertion follows.

The sequence $(\lambda _n)_n$ is defined in terms of $(m_n)_n$ from
(\ref{eq:fix}) by
\begin{equation}\label{eq:sl}
\lambda _0=0, \quad \lambda _{n+1}=1/m_n,\quad n\ge 0.
\end{equation}
By (\ref{eq:Mellin}) and (\ref{eq:fF=1}) we clearly have 
\begin{equation}\label{eq:m-la}
m_n=F(n),\;\lambda_n=f(n),\quad n\ge 0,
\end{equation}
hence by (\ref{eq:it})
\begin{equation}\label{eq:fes1}
\lambda_{n}=\psi(\lambda_{n+1}),\quad n\ge 0,
\end{equation}
which can be reformulated to
\begin{equation}\label{eq:fes2}
 \lambda _{n+1}=\frac{1}{2}\left( \lambda _{n}+\sqrt {\lambda _n^2+4}\right),
\quad n\ge 0.
\end{equation}

The following result is easy and the proof is left to the reader.

\begin{lem}\label{thm:Y} Defining 
\begin{equation}\label{eq:Y} 
Y_n=(\psi^{\circ n})^{-1}(\{0\})=\{z\in\mathbb C\mid \psi^{\circ
  n}(z)=0\},
\end{equation}
i.e.
$$
Y_0=\{0\},\quad Y_1=\{-1,1\},\quad Y_2=\{(\pm 1\pm\sqrt{5})/2\},\ldots
$$
we have for $n\ge 1$ 
\begin{enumerate}
\item[(i)] $\psi(Y_{n})=Y_{n-1},\;Y_{n}=\psi^{\circ -1}(Y_{n-1})$,
\item[(ii)] The set of poles of $\psi^{\circ n}$ is
  $\cup_{j=0}^{n-1}Y_j$,
\item[(iii)] $Y_n$ consists of $2^n$ real numbers and is symmetric
  with respect to zero.
\item[(iv)]  The function $\psi^{\circ n}$
  is strictly increasing from $-\infty$ to $\infty$ in each of the
  $2^n$ intervals in which $\cup_{j=0}^{n-1}Y_j$ divides $\mathbb
  R$. There is exactly one zero of $\psi^{\circ n}$ in each of these
 intervals, and these
  zeros form the set $Y_n$. 
\end{enumerate}
\end{lem}

We write $Y_n=\{ \alpha _{n,k}: k=1,\ldots , 2^n\}$ arranged in
increasing order ($n\ge 1$):
$$
\alpha_{n,1}<\alpha_{n,2}<\cdots<\alpha_{n,2^{n-1}}<0<\alpha_{n,2^{n-1}+1}
<\cdots<\alpha_{n,2^n}.
$$
It is easy to see that
 $-\alpha_{n,1}=\alpha_{n,2^n}=\lambda_n$ for $n\ge 0$.

\begin{prop}\label{thm:unionY} The set
$$
\cup_{p=0}^\infty Y_p=\left\{\alpha_{p,k} \mid p\ge
  0,k=1,\ldots,2^p\right\}
$$
is dense in $\mathbb R$.
\end{prop}

\begin{proof} The set in question is the so-called backward orbit of 0
for $\psi$, and since $0\in J(\psi)$ the result follows
 by \cite[Theorem 4.2.7]{bear}.
\end{proof}

We next give some asymptotic properties of the sequence $(\lambda
_n)_n$ and the function $f$:

\begin{lem} \label{thm:aps}
\begin{enumerate}
\item $\displaystyle \sqrt n\le \lambda _n\le \sqrt{2n}$, $n\ge 0$.
\item $(\lambda_n)_n$ is an increasing divergent sequence
and $\lambda_{n+1}/\lambda_n$ is decreasing with
$\displaystyle \lim _{n\to\infty} \frac{\lambda _{n+1}}{\lambda _n}=1$.
\item $\displaystyle \lim _{n\to\infty} (\lambda ^2_{n+1}-\lambda ^2_n)=2$.
\item $\displaystyle \lim _{n\to\infty} \frac{\lambda ^2_n}{n}=2$.
\item $\displaystyle \lim _{n\to\infty} \frac{\lambda _n^2-2n}{\log
    n}=-\frac{1}{2}$.
\item $\lim _{s\to \infty}f(s)/\sqrt{2s}=1$.
\item $\lim _{s\to \infty}f'(s)\sqrt{2s}=1$.
\end{enumerate}
\end{lem}

\begin{proof}

1. These inequalities follow easily from (\ref{eq:fes2}) using
induction on $n$.

2. The sequence $(\lambda _n)_n$  increases to infinity since it is the
reciprocal of the Hausdorff moment sequence $(m_n)_n$. By the
Cauchy-Schwarz inequality $m_n^2\leq m_{n-1}m_{n+1}$, which proves
that $(\lambda_{n+1}/\lambda_n)_n$ is decreasing.
The limit follows now easily
 from (\ref{eq:fes2}).

3. Using (\ref{eq:fes1}) we can write
$$
\lambda _{n+1}^2-\lambda _n^2=\frac{\lambda _{n+1}+\lambda _n}{\lambda_{n+1}}=
1+\frac{\lambda _n}{\lambda _{n+1}},
$$
and it suffices to apply part 2.

4. is a consequence of part 3 and the following version of the Stolz
 criterion going back to \cite{Stolz}:

\begin{lem}\label{thm:Stolz}
Let $(a_n)_n,(b_n)_n$ be real sequences, where $(b_n)_n$ is strictly
increasing  tending to infinity. Then
$$
\lim_{n\to\infty} \frac{a_{n+1}-a_n}{b_{n+1}-b_n}=L\Rightarrow
 \lim_{n\to\infty} \frac{a_n}{b_n}=L.
$$
\end{lem}

5. follows  by using again the Stolz criterion  and taking into account that
$$
\frac{\lambda _{n+1}^2-\lambda _n^2-2}{\log \frac{n+1}{n}}=
\frac{\lambda _{n+1}^2-\lambda _n^2-2\lambda_{n+1}^2+2\lambda
  _{n+1}\lambda_n}{\log \frac{n+1}{n}}
$$
$$
=-\frac{(\lambda _{n+1}-\lambda _n)^2}{\log \frac{n+1}{n}}
 =-\frac{1}{n\log \frac{n+1}{n}}\frac{n}{\lambda _{n+1}^2
 }\to -\tfrac12.
$$

6. Since $f$ is increasing and $f(n)=\lambda _n$, the assertion follows
from part 4.

7. We write $f(n+1)-f(n)=f'(t_n)$, for a certain
 $t_n\in (n,n+1)$.
Since $f'$ is decreasing ($f'(s)$ is completely monotonic),
 part 7 follows if we prove that $f'(t_n)\sqrt {2t_n}$
 tends to $1$ as $n$ tends to $\infty $.
However, using the recursion formula for $(\lambda _n)_n$, we get
$$
f'(t_n)\sqrt {2t_n}=(\lambda _{n+1}-\lambda_n)\sqrt {2t_n}=
\frac{\sqrt {2(n+1)}}{\lambda _{n+1}}
\frac{\sqrt {2t_n}}{\sqrt {2(n+1)}},
$$
and it suffices to apply part 4.
\end{proof}

{\sl Proof of Theorem \ref{thm:main-b}}.

\medskip
We have already proved the properties (i) and (iii). To see (ii) we
notice that $f=\mathcal B(\mu)$ is a
Bernstein function, and therefore $1/f$ is completely monotonic.
Every completely monotonic function is logarithmically convex. For
these statements see e.g. \cite[\S 14]{bf}.

Suppose next that $\tilde f$ is a function satisfying (i)-(iii). Since
$\tilde f(1)=1=\lambda_1$,  we see by (iii) and (\ref{eq:fes1})
 that $\tilde f(n)=\lambda_n$ for $n\ge
1$. Equation (\ref{eq:iterate}) is equivalent with
 \begin{equation}\label{eq:iterate1}
\tilde f(s)=\lim_{n\to\infty}\psi^{\circ n}\left(\lambda_n\left(
\frac{\lambda_{n+1}}{\lambda_n}\right)^s\right),
\end{equation}  
and if we prove this equation for $0<s\le 1$, then $\tilde f$ is uniquely
 determined on $]0,1]$ and hence by (iii) for all $s>0$.

We prove that the limit in (\ref{eq:iterate1}) exists and
coincides with
 $\tilde f(s)$ for $0<s\le 1$. This is clear for $s=1$ since
$\psi^{\circ n}(\lambda_{n+1})=1$ for $n\geq 0$.

For any convex function $\phi $ on $]0,\infty[$  we have for $0< s\leq 1$
and $n\ge 2$ 
$$
\phi (n)-\phi (n-1)\le \frac{\phi (n+s)-\phi (n)}{s}\le \phi (n+1)-\phi (n).
$$
By taking $\phi =\log (1/\tilde f)$, which it is convex by assumption, we get
$$
\log \frac{\lambda _{n-1}}{\lambda _n}\le \frac{1}{s}\log \frac{\tilde
  f(n)}{\tilde f(n+s)}\le \log \frac{\lambda _{n}}{\lambda _{n+1}};
$$
that is
$$
\left(\frac{\lambda _{n-1}}{\lambda _n}\right) ^s\le \frac{\lambda
  _n}{\tilde f(n+s)}\le
\left( \frac{\lambda _{n}}{\lambda _{n+1}}\right) ^s,
$$
which finally gives:
$$
\lambda _n\left( \frac{\lambda _{n+1}}{\lambda _{n}}\right) ^s
\le \tilde f(n+s)\le
\lambda _n \left(\frac{\lambda _{n}}{\lambda _{n-1}}\right) ^s, \quad 0<s<1.
$$

Using that $\psi$ is increasing
on $]0,\infty[$, we get by applying $\psi^{\circ n}$  to the previous
 inequality
$$
\psi ^{\circ n} (b_n(s))\le \tilde f(s)=\psi^{\circ n}(\tilde f(n+s))\le \psi ^{\circ n} (a_n(s)),
$$
where we have introduced
$$
a_n(s)=\lambda _n \left(\frac{\lambda _{n}}{\lambda _{n-1}}\right) ^s,\quad
b_n(s)=\lambda _n\left( \frac{\lambda _{n+1}}{\lambda _{n}}\right) ^s.
$$

It is now enough to prove that
$$
\lim_{n\to\infty} (\psi ^{\circ n}
 (a_n(s))-\psi ^{\circ n} (b_n(s))=0.
$$
By applying the mean value theorem, we get for a  certain
$w\in ]b_n(s),a_n(s)[$ that
$$
\psi ^{\circ n} (a_n(s))-\psi ^{\circ n} (b_n(s))=
(a_n(s)-b_n(s))(\psi ^{\circ n})' (w)
$$
$$
=(a_n(s)-b_n(s))\psi '(\psi ^{\circ n-1} (w))\psi '(\psi ^{\circ n-2} (w))
\cdots \psi '(w).
$$

Since $\lambda _n<b_n(s)<w<a_n(s)$, we get
$\lambda _{n-k}<\psi ^{\circ k}(b_n(s))<\psi ^{\circ k}(w)$, $k=0,1,
\ldots ,n$, hence
\begin{eqnarray*}
\lefteqn{\vert \psi ^{\circ n} (a_n(s))-\psi ^{\circ n} (b_n(s)\vert}\\
 &\le &
 \vert a_n(s)-b_n(s)\vert \prod _{k=0}^{n-1}\vert \psi '(\psi ^{\circ k} (w))\vert\\
&\le &
\vert a_n(s)-b_n(s)\vert \prod _{k=0}^{n-1}\left( 1+\frac{1}{\lambda _{n-k}^2}\right)\\
&=&
\lambda _n\left( \left(\frac{\lambda _{n}}{\lambda _{n-1}}\right)^s-
\left( \frac{\lambda _{n+1}}{\lambda _{n}}\right) ^s\right) \prod _{k=1}^{n}\left( 1+\frac{1}{\lambda _{k}^2}\right)\\
&\le& \lambda _n\left( \left(\frac{\lambda _{n}}{\lambda_{n-1}}
\right)^s-\left( \frac{\lambda _{n+1}}{\lambda _{n}}\right) ^s\right)
 \prod _{k=1}^{n}\left( 1+\frac{1}{k}\right)\\
&=& (n+1)\lambda _n\left( \left(\frac{\lambda _{n}}{\lambda _{n-1}}\right)^s-
\left( \frac{\lambda _{n+1}}{\lambda _{n}}\right) ^s\right),
\end{eqnarray*}
where we have used $\sqrt{k}\leq\lambda_k$ from Lemma \ref{thm:aps}
part 1.

Using that $(x^s-y^s)\le s(x-y)$ for $1<y<x$ and $0<s\leq 1$, we get
$$
\vert \psi ^{\circ n} (a_n(s))-\psi ^{\circ n} (b_n(s))\vert
\le s(n+1)\lambda _n\left( \frac{\lambda _{n}}{\lambda _{n-1}}-
 \frac{\lambda _{n+1}}{\lambda _{n}}\right),
$$
and by (\ref{eq:fes2}) we finally get
\begin{eqnarray*}
\lefteqn{\vert \psi ^{\circ n} (a_n(s))-\psi ^{\circ n} (b_n(s)\vert}\\
&\le&
  \tfrac12 s(n+1)\lambda _n\left(
\left( 1+\sqrt {1+ \frac{4}{\lambda _{n-1}^2}}\right)-
\left( 1+\sqrt {1+ \frac{4}{\lambda _{n}^2}}\right) \right)\\
&=&
\tfrac12 s(n+1)\lambda _n\left( \sqrt {1+ \frac{4}{\lambda _{n-1}^2}}-
\sqrt {1+ \frac{4}{\lambda _{n}^2}}\right)\\
&=&\frac{2s(n+1)\lambda _n\displaystyle \left(\frac{1}{\lambda _{n-1}^2}-
\frac{1}{\lambda _{n}^2}\right)}{\sqrt {1+ \frac{4}{\lambda _{n-1}^2}}+
\sqrt {1+ \frac{4}{\lambda _{n}^2}}}
\le \frac{s(n+1)}{\lambda _n\lambda ^2_{n-1}}(\lambda _n^2-\lambda ^2_{n-1}),
\end{eqnarray*}
which tends to zero by part 2, 3 and 4 of Lemma \ref{thm:aps}.
\hfill$\square$

\medskip

For each real number $s$, we define the sequence $(\lambda _n (s))_n$ by $\lambda _0(s)=s$ and
\begin{equation}\label{eq:lambdas}
\lambda _{n+1}(s)=\frac{\lambda _n(s)+\sqrt{\lambda _n(s)^2+4}}{2},
 \quad n\ge 0.
\end{equation}
Notice that $\lambda _{n+1}(s)$ is the positive root of
 $z^2-\lambda _n(s)z-1=0$ and that
\begin{equation}\label{eq:psilambda}
\psi(\lambda_{n+1}(s))=\lambda_n(s).
\end{equation} 
Therefore,
if $s\in Y_l$ then $\lambda _n(s)\in Y_{l+n}$, and for $s=0$ we
 have $\lambda _n(0)=\lambda _n$, $n\ge 0$. Furthermore,
 $\lambda_n(\lambda_l(s))=\lambda_{n+l}(s).$

\begin{defn}\label{thm:def}
For integers $k,l\ge 0$ we denote by $r(k,l)$  the unique solution
$x\in\{1,2,\ldots,2^l\}$ of the congruence equation $x\equiv k \mod 2^l$.
\end{defn}

\begin{lem}\label{thm:step1} For $p\ge 1, k=1,2,\ldots,2^p$ we have
\begin{enumerate}
\item[(i)] $\psi(\alpha_{p,k})=\alpha_{p-1,r(k,p-1)}$.
\item[(ii)] $\psi^{\circ l}(\alpha_{p,k})= \alpha_{p-l,r(k,p-l)}$
  for $l=0,1,\ldots,p$.
\end{enumerate}
\end{lem}

\begin{proof} Since $\psi(Y_p)=Y_{p-1}$ and $\psi$ is strictly increasing
  mapping $\left]-\infty,0\right[$ onto $\mathbb R$, we see that
$$
\psi(\alpha_{p,k})=\alpha_{p-1,k},\quad k=1,2,\ldots,2^{p-1},
$$
and since similarly $\psi$ maps $\left]0,\infty\right[$ onto
$\mathbb R$ we get
$$
\psi(\alpha_{p,k})=\alpha_{p-1,j},\quad
k=2^{p-1}+j,\;j=1,2,\ldots,2^{p-1}.
$$
In the first case $k=r(k,p-1)$ and in the second case $j=r(k,p-1)$ so
the assertion (i) follows.

The assertion (ii) is clear for $l=0$ and $l=p$ and follows for $l=1$
by (i). Assuming (ii) for some $l$ such that $1\le l\le p-2$ we get by
(i)
$$
\psi^{\circ(l+1)}(\alpha_{p,k})=\psi(\alpha_{p-l,r(k,p-l)})
=\alpha_{p-l-1,j},
$$
where $j:=r(r(k,p-l),p-l-1)$. By definition
$$
k\equiv r(k,p-l) \mod 2^{p-l},\;1\le r(k,p-l)\le 2^{p-l}
$$
$$
j\equiv r(k,p-l) \mod 2^{p-l-1},\; 1\le j\le 2^{p-l-1}.
$$
The first congruence also holds $\mod 2^{p-l-1}$, hence $j\equiv k
\mod 2^{p-l-1}$ and finally $j=r(k,p-l-1)$.
\end{proof}

\begin{cor}\label{thm:step2} For a zero $\xi_{p,k}$ of $f$ we have
\begin{enumerate}
\item[(i)] $f(\xi_{p,k}+l)=\alpha_{l,r(k,l)},\; l=0,1,\ldots,p,$
\item[(ii)] $f(\xi_{p,k}+l)=\lambda_{l-p}(\alpha_{p,k}),\;l=p+1,p+2,
\ldots$,
where $\lambda_n(s)$ is defined in (\ref{eq:lambdas}).
\end{enumerate}
\end{cor}

\begin{proof} We first prove (i) for $l=p$, i.e. that
$f(\xi_{p,k}+p)=\alpha_{p,k}$ since $r(k,p)=k$. Note that by
 (\ref{eq:it}) we have
 $$
 \psi^{\circ p}(f(\xi_{p,k}+p))=f(\xi_{p,k})=0,
 $$
hence $f(\xi_{p,k}+p)\in Y_p$. On the other hand
  $\xi_{p,k}+p\in\left]-1,0\right[$, and since $f$ is strictly
  increasing satisfying
  $f(]-1,0[)=]-\infty,0[$, we see that $f(\xi_{p,k}+p),k=1,2,\ldots, 2^{p-1}$
  describe $2^{p-1}$ negative numbers in $Y_p$ in increasing order. Therefore,
  $f(\xi_{p,k}+p)=\alpha_{p,k}, k=1,2,\ldots,2^{p-1}$.

  By Lemma \ref{thm:step1} and (\ref{eq:it}) we then get for $0\le l\le p$
  $$
  f(\xi_{p,k}+l)=\psi^{\circ(p-l)}(f(\xi_{p,k}+p))=\psi^{\circ(p-l)}
  (\alpha_{p,k})=\alpha_{l,r(k,l)}.
  $$
  Clearly $0<f(\xi_{p,k}+p+1)\in Y_{p+1}$ and $\alpha_{p,k}=
  \psi(f(\xi_{p,k}+p+1))$, hence $f(\xi_{p,k}+p+1)=\lambda_1(\alpha_{p,k})$
  by definition of $\lambda_1(s)$. The assertion (ii) follows easily by induction.
\end{proof}

\begin{thm}\label{thm:xi-rho}
The numbers
$\xi _{p,k}$, $\rho _{p,k}$, $p\ge 1$, $k=1,\ldots ,2^{p-1}$ and
$\rho_0$  from Theorem \ref{thm:main-d} are given by the following formulas:
\begin{equation} \label{eq:xi}
\xi _{p,k}=\lim _{N\to \infty} \sqrt{2N}\left(
\sum _{l=1}^{p}\frac{1}{\alpha_{l,r(k,l)}}+
\sum _{l=1}^{N-p}\frac{1}{\lambda _l(\alpha _{p,k})} -\lambda _N\right),
\end{equation}

\begin{equation}\label{eq:rho}
\rho _{p,k}=\prod _{l=1}^{p}\left( 1+
\frac{1}{\alpha _{l,r(k,l)}^2}\right)^{-1}
\lim _{N\to \infty} \sqrt {2N}\prod _{l=1}^{N}\left( 1+
\frac{1}{\lambda _l^2(\alpha _{p,k})}\right)^{-1},
\end{equation}

\begin{equation}\label{eq:rho-0}
\rho _0=
\lim _{N\to \infty}\sqrt {2N}\prod _{l=1}^{N}\left(
    1+\frac{1}{\lambda _l^2}\right)^{-1}.
\end{equation}
\end{thm}

\begin{proof} By applying $N$ times the functional equation
 (\ref{eq:z,z+1}) for the function
 $f$ and using Corollary \ref{thm:step2} , we have for $p<N$:
$$
0=f(\xi _{p,k})=f(\xi _{p,k}+N)-\sum _{l=1}^N\frac{1}{f(\xi _{p,k}+l)}
$$
$$
=
f(\xi _{p,k}+N)-\left( \sum _{l=1}^{p}\frac{1}{\alpha _{l,r(k,l)}}+
\sum _{l=1}^{N-p}\frac{1}{\lambda _l(\alpha _{p,k})}\right).
$$
Writing
$$
y_{N,p,k}=\sum _{l=1}^{p}\frac{1}{\alpha _{l,r(k,l)}}+
\sum _{l=1}^{N-p}\frac{1}{\lambda _l(\alpha _{p,k})},
$$
we get $f(\xi_{p,k}+N)=y_{N,p,k}$. For $N\to\infty$ it follows by part 6 of
Lemma \ref{thm:aps}  that
$y_{N,p,k}\sim \sqrt {2N}$. Since $f$ is a strictly increasing 
bijection of $(-1,+\infty )$ onto $\mathbb R$, we can consider
its inverse $f^{-1}$. Then we have $N=f^{-1}(\lambda _N)$, hence
$\xi _{p,k}=f^{-1}(y_{N,p,k})-f^{-1}(\lambda _N)$. Since $\xi _{p,k}$
 is negative and $f$ is
increasing, we deduce that $y_{N,p,k}<\lambda _N$. This gives for a certain
number $\sigma _{N,p,k}\in \left]y _{N,p,k}, \lambda _N\right[$ that
$$
\xi _{p,k}=f^{-1}(y_{N,p,k})-f^{-1}(\lambda_N)=(f^{-1})'(\sigma
_{N,p,k})
(y_{N,p,k}-\lambda_N)=
\frac{y_{N,p,k}-\lambda_N}{f'(\eta  _{N,p,k})},
$$
where we have written $\eta _{N,p,k}=f^{-1}(\sigma _{N,p,k})$.
Clearly $\eta _{N,p,k}\in \left]\xi_{p,k}+N, N\right[$.

Taking into account that $\lim _{s\to \infty}f'(s)\sqrt {2s}=1$ 
(part 7 of Lemma \ref{thm:aps}),
we have
$$
\xi _{p,k}=\lim _N \sqrt{2N}\left(y_{N,p,k}-\lambda_N\right),
$$
that is, (\ref{eq:xi}) holds.

 The number $f'(\xi _{p,k})$ can be computed as follows:
Deriving the functional equation (\ref{eq:z,z+1}) for $f$, we get
$$
f'(z)=f'(z+1)\left(1+\frac{1}{f^2(z+1)}\right)
$$
hence by iteration
\begin{equation}\label{eq:f'ite}
f'(z)=f'(z+N)\prod _{l=1}^N\left( 1+\frac{1}{f^2
(z+l)}\right).
\end{equation}
Using Corollary \ref{thm:step2} and $\lim _{s\to \infty}f'(s)\sqrt {2s}=1$,
(Lemma \ref{thm:aps}, part 7) we get for $z=\xi_{p,k}$
$$
f'(\xi _{p,k}) =\prod _{l=1}^{p}\left( 1+\frac{1}{\alpha _{l,r(k,l)}^2}\right)
\lim _{N\to \infty} \frac{1}{\sqrt {2N}}\prod _{l=1}^{N}\left(
  1+\frac{1}{\lambda _l^2(\alpha _{p,k})}\right) ,
$$
and since $\rho_{p,k}=1/f'(\xi_{p,k})$ by (\ref{eq:rho-f}), we see that
(\ref{eq:rho}) holds.

Applying (\ref{eq:f'ite}) for $z=0$, we get
$$
f'(0)=f'(N)\prod_{l=1}^N\left(1+\frac{1}{\lambda_l^2}\right),
$$
and (\ref{eq:rho-0}) follows by (\ref{eq:rho-f}) and
 $\lim_{N\to\infty}f'(N)\sqrt{2N}=1$.
\end{proof}
We give some values of the numbers of Theorem \ref{thm:xi-rho}:  
$$
\boxmatrix \\
& \rho_0=0.68\ldots && \xi_0=0 &\\
& \rho_{1,1}=0.14\ldots && \xi_{1,1}=-1.46\ldots &\\
& \rho_{2,1}=0.06\ldots && \xi_{2,1}=-2.61\ldots &\\
& \rho_{2,2}=0.05\ldots && \xi_{2,2}=-2.33\ldots &\\ 
\endboxmatrix
$$

\begin{thm}\label{thm:one} The density $\mathcal D$ given by
(\ref{eq:D}) satisfies
$$
\mathcal D(t)\sim \frac{1}{\sqrt{2\pi(1-t)}}\textrm{ for } t\to 1.
$$
\end{thm}

\begin{proof} By formula (\ref{eq:fF=1}) and Lemma \ref{thm:aps} part
  6 we get
$$
F(s)=\int_0^1 t^s\mathcal D(t)\,dt \sim \frac{1}{\sqrt{2s}},\quad
s\to\infty,
$$
or
$$\int_0^\infty e^{-us}\mathcal D(e^{-u})e^{-u}\,du \sim
\frac{1}{\sqrt{2s}},
\quad s\to\infty.
$$

By the Karamata Tauberian theorem, cf. \cite[Theorem 1.7.1$'$]{bgt},  we
get
$$
\int_0^t \mathcal D(e^{-u})e^{-u}\,du \sim \sqrt{\frac{2t}{\pi}},\quad
t\to 0,
$$
and since $\mathcal D$ is increasing we can use the Monotone Density
theorem, cf. \cite[Theorem 1.7.2b]{bgt}, to conclude that
$$
\mathcal D(e^{-u})e^{-u} \sim\frac{1}{\sqrt{2\pi u}},\quad u\to 0,
$$
which is equivalent to the assertion.
\end{proof}

\section{Miscellaneous about the fixed point}

The fixed point sequence $(m_n)_n$ given by (\ref{eq:fix}) satisfies
$m_{n+1}=\Phi(m_n)$ with
$$
\Phi(x)=\frac{\sqrt{4x^2+1}-1}{2x},\quad x>0.
$$
This makes it possible to express $(m_n)_n$ as iterates of $\Phi$, viz.
$$
m_n=\Phi^{\circ n}(1).
$$
From Lemma \ref{thm:aps} part 4 we get the asymptotic behaviour of
$m_n$ as
$$
m_n\sim\frac{1}{\sqrt{2n}},\quad n\to\infty.
$$
This behaviour can also be deduced from a general result about
iteration, cf. \cite[p.175]{db}. The authors want to thank Bruce
Reznick for this reference as well as the following description of
$(m_n)_n$.

\begin{prop}\label{thm:Rez} Define $h_n\in\left]0,\pi/4\right]$ by
  $\tan h_n=m_n$ and let
$$
G(x)=\frac12\arctan(2\tan x),\quad |x|<\frac{\pi}2.
$$
Then
$$
h_n=G^{\circ n}(\frac{\pi}{4}).
$$
\end{prop}

\begin{proof} We have
$$
\tan h_n=m_n=\frac{m_{n+1}}{1-m_{n+1}^2}=\frac{\tan h_{n+1}}{1-\tan^2
  h_{n+1}}=\frac12\tan(2h_{n+1}),
$$
hence $h_{n+1}=G(h_n)$ and the assertion follows.
\end{proof}

A Hausdorff moment sequence $(a_n)_n$ is called {\sl infinitely divisible}
if $(a_n^{\alpha})_n$ is a Hausdorff moment sequence for all $\alpha>0$. If
$a_n=\int_0^1 t^n\,d\nu(t),n\ge 0$ then $(a_n)_n$ is infinitely divisible
 if and only if  $\nu$ is infinitely divisible for the product
 convolution $\tau\diamond\nu$ of measures $[0,\infty[$ defined by
$$
\int g\,d\tau\diamond\nu=\int\int g(st)\,d\tau(s)\,d\nu(t).
$$
For a general study of these concepts
   see \cite{ty},\cite{b2},\cite{b3}. In case the measure $\nu$ does
   not charge $0$, the notion is the classical infinite divisibility on
   the locally compact group $]0,\infty[$ under multiplication.
                                     
\begin{prop}\label{thm:infdiv} Hausdorff moment sequences of the form 
(\ref{eq:tms}) are infinitely divisible.
\end{prop}

\begin{proof} Let $\nu\ne 0$ be a positive measure on $[0,1]$ and let
$a_n=\int t^n\,d\nu(t),n\ge 0$ be the corresponding Hausdorff moment sequence.
Let $\alpha>0$ be fixed. We shall prove that
 $((a_0+a_1+\cdots+a_n)^{-\alpha})_n$
is a Hausdorff moment sequence.

 For $0<c<1$ we denote by $\nu_c=\nu|[0,c[+\nu(\{1\})\delta_c$, where
 the first term denotes the restriction of $\nu$ to $[0,c[$. Then
$\lim_{c\to 1}\nu_c=\nu$ weakly and in particular for each $n\ge 0$
$$
a_n(c):=\int_0^1 t^n\,d\nu_c(t)\to a_n\;\;\mbox{for}\;\; c\to 1.
$$
It therefore suffices to prove that
\begin{equation}\label{eq:anc}
((a_0(c)+a_1(c)+\cdots+a_n(c))^{-\alpha})_n
\end{equation}
is a Hausdorff moment sequence. 
 By a simple calculation we find
\begin{eqnarray*}
\left(\sum_{k=0}^n
  a_k(c)\right)^{-\alpha}&=&\left(\int_0^1\frac{1-t^{n+1}}{1-t}\,d\nu_c(t)
\right)^{-\alpha}\\
&=&
\left(\int_0^1\frac{d\nu_c(t)}{1-t}-\int_0^1 t^n\frac{t\,d\nu_c(t)}{1-t}
\right)^{-\alpha}=H(\tau_n),
\end{eqnarray*}
where
$$
\tau_n=\int_0^1 t^n\frac{t\,d\nu_c(t)}{1-t},\quad
H(z)=\left(\int_0^1\frac{d\nu_c(t)}{1-t} -z\right)^{-\alpha}.
$$
The function $H$ is clearly holomorphic in
$$
|z|<\int_0^1\frac{d\nu_c(t)}{1-t}
$$
with non-negative coefficients in the power series. Applying Lemma 2.1 in
\cite{bd2}, shows that (\ref{eq:anc}) is a Hausdorff moment sequence.
\end{proof}

\begin{cor}\label{thm:infdivfix} The fixed point sequence $(m_n)_n$ is
  infinitely divisible.
\end{cor}

\begin{rem}\label{thm:infdiv1} {\rm By Corollary \ref{thm:infdivfix}
the fixed point measure $\mu$ is infinitely divisible  for the
 product convolution. The image measure $\eta$ of $\mu$ under
 $\log(1/t)$ is an infinitely divisible probability measure in the
 ordinary sense, because $\log(1/t)$ maps products to sums. The
 measure $\eta$ has the density
\begin{equation}\label{eq:D(e-u)}
\mathcal D(e^{-u})e^{-u}=\rho_0
e^{-u}+\sum_{p=1}^\infty\sum_{k=1}^{2^{p-1}}
\rho_{p,k}e^{-u(1-\xi_{p,k})},\quad u>0
\end{equation}
with respect to Lebesgue measure on the half-line. Since
(\ref{eq:D(e-u)}) is clearly a completely monotonic density, the
infinite divisibility of $\eta$ is also a consequence of the
Goldie-Steutel theorem, see \cite[Theorem 10.7]{S:vH}. These remarks
also show that Corollary \ref{thm:infdivfix} can be inferred from the complete
monotonicity of (\ref{eq:D(e-u)}) via the Goldie-Steutel theorem. The
formula
$$
\int_0^\infty e^{-us}\,d\eta(u)=\int_0^1 t^s\,d\mu(t)=F(s)=e^{-\log
  f(s+1)},\quad s\ge 0
$$
shows that $\log f(s+1)$ is the Bernstein function associated with the
 convolution semigroup $(\eta_t)_{t>0}$ of probability measures on
 the half-line such that $\eta_1=\eta$, see \cite[p. 68]{bf}.}
\end{rem}

\begin{rem}\label{thm:infdiv2} {\rm Let $\mathcal H_I$ denote the set
of normalized infinitely divisible Hausdorff moment sequences. By 
Proposition \ref{thm:infdiv} we have $T(\mathcal H)\subseteq \mathcal
H_I$. We claim that this inclusion is proper. In fact, it is easy to
see that $T:\mathcal H\to T(\mathcal H)$ is one-to-one, and that
$$
T^{-1}(\mathbf{b})_n=\frac{1}{b_n}-\frac{1}{b_{n-1}},\quad n\ge 1,
$$
for $\mathbf{b}=(b_n)_n\in T(\mathcal H)$. It follows that
$$
T(\mathcal H)=\{\mathbf{b}\in\mathcal H \mid \left(\frac{1}{b_n}-\frac{1}{b_{n-1}}
\right)_n\in\mathcal H\}.
$$
(Here $1/b_n-1/b_{n-1}=1$ for $n=0$.) Then $\mathbf{b}\in\mathcal
H_I\setminus T(\mathcal H)$ if we define $b_n=1/(n+1)^2$.}
\end{rem}
 
The functions $f,F$ being holomorphic in $\Re z>-1$ with a pole at
$z=-1$, they have power series expansions
\begin{equation}\label{eq:power}
F(z)=1+\sum_{n=1}^\infty a_nz^n,\quad f(z)=\sum_{n=1}^\infty
b_nz^n,\quad |z|<1,
\end{equation}
and the radius of convergence is 1 for both series.

\begin{prop}\label{thm:power} The coefficients in (\ref{eq:power}) are
  given for $n\ge 1$ by
\begin{eqnarray*}
a_n&=&\frac{1}{n!}\int_0^1(\log t)^n\,d\mu(t)=(-1)^n\left(\rho_0+
\sum_{p=1}^\infty\sum_{k=1}^{2^{p-1}}\frac{\rho_{p,k}}{(1-\xi_{p,k})^{n+1}}\right),\\
b_n&=&-\frac{1}{n!}\int_0^1\frac{(\log
  t)^n}{1-t}\,d\mu(t)\\
&=&(-1)^{n-1}\left(\rho_0\zeta(n+1,0)+
\sum_{p=1}^\infty\sum_{k=1}^{2^{p-1}}\rho_{p,k}\zeta(n+1,-\xi_{p,k})\right),
\end{eqnarray*}
where 
$$
\zeta(s,a)=\sum_{n=1}^\infty\frac{1}{(n+a)^s},\quad s>1,a>-1
$$
is the Hurwitz zeta function.
\end{prop}

\begin{proof} The formula for $a_n$ follows from (\ref{eq:Mellin}) and
(\ref{eq:F}), and the formula for $b_n$ follows from (\ref{eq:Bern}) and
(\ref{eq:f}).
\end{proof}

\medskip
{\bf Acknowledgement.} The authors wish to thank Henrik L. Pedersen for
help with producing the graph of $f$.

\end{document}